\documentclass[final, letter]{elsart1p}
\usepackage{setspace}
\usepackage{ifpdf}
\usepackage{graphicx,amssymb,lineno}
\ifpdf
\usepackage[%
  pdftitle={Instructions for use of the document class
    elsart},%
  pdfauthor={Simon Pepping},%
  pdfsubject={The preprint document class elsart},%
  pdfkeywords={instructions for use, elsart, document class},%
  pdfstartview=FitH,%
  bookmarks=true,%
  bookmarksopen=true,%
  breaklinks=true,%
  colorlinks=true,%
  linkcolor=blue,anchorcolor=blue,%
  citecolor=blue,filecolor=blue,%
  menucolor=blue,pagecolor=blue,%
  urlcolor=blue]{hyperref}
\else
\usepackage[%
  breaklinks=true,%
  colorlinks=true,%
  linkcolor=blue,anchorcolor=blue,%
  citecolor=blue,filecolor=blue,%
  menucolor=blue,pagecolor=blue,%
  urlcolor=blue]{hyperref}
\fi

\makeatletter
\def\elsartstyle{%
    \def\normalsize{\@setfontsize\normalsize\@xiipt{14.5}}
    \def\small{\@setfontsize\small\@xipt{13.6}}
    \let\footnotesize=\small
    \def\large{\@setfontsize\large\@xivpt{18}}
    \def\Large{\@setfontsize\Large\@xviipt{22}}
    \skip\@mpfootins = 18\p@ \@plus 2\p@
    \normalsize
}
\@ifundefined{square}{}{}
\makeatother

\newtheorem{theorem}[section]{\bf Theorem}
\newtheorem{definition}{\bf Definition}
\newtheorem{corollary}{\bf Corollary}
\newtheorem{exam}{\bf Example}

\newtheorem{lemma}{\bf Lemma}
\newcommand{\dm}[1]{ {\displaystyle{#1} } }
\def \R{{\mathbb R}}
\def \C{{\mathbb C}}

\def \E{{\bf E}}

\def \rank{{\bf rank}}

\newcommand{\norm}[1]{\|{#1}\|}

\def \sign{\mathrm{sign}}

\def \pf{{\bf Proof: }}

\def \lam{\lambda}

\def \Lam{\Lambda}

\def \pf{{\bf Proof: }}

\def \cnn{\C^{n\times n}}

\def \norm{\|\cdot\|}
\begin{document}
\begin{frontmatter}
\title{\bf The Stochastic Logarithmic Norm for Stability Analysis of
	Stochastic Differential Equations} 
\author[a,b]{Sk. Safique Ahmad} \ead{safique@gmail.com}
\and
\author[a]{Nagalinga Rajan} \ead{rajan@rishi.serc.iisc.ernet.in}
\and
 \author[a,b]{Soumyendu Raha}  \ead{raha@serc.iisc.ernet.in}
 \thanks[a]{Scientific Computation Laboratory, Supercomputer Education and Research Centre, Indian Institute of Science, Bangalore 560012, India. \\ The last author is the corresponding author.}
\thanks[b]{This work was partially supported by the Naval Research Board, DRDO, Government of India with grant number DNRD/05/4003/NRB/88.}
\date{\today}
\begin{abstract}
To analyze the stability of It\^o stochastic differential equations with multiplicative noise, 
we  introduce the stochastic logarithmic norm. 
The logarithmic norm was originally introduced by G. Dahlquist in 1958 as a tool to study the growth of solutions to ordinary differential equations and for estimating the error growth in discretization methods for their approximate solutions. We extend the concept to the stability analysis of It\^o stochastic differential equations with multiplicative noise. Stability estimates for linear It\^o SDEs using the one, two and $\infty$-norms in the
$l$-th mean, where $1 \leq l < \infty $, are derived and the application of the stochastic logarithmic norm is illustrated with examples.
\end{abstract}
%\von\noin {\bf AMS subject classification:}.........
\end{frontmatter}

\section{Introduction}
For $A\in \C^{n\times n}$ and $X(t)\in \C^n,$ we consider the
ordinary differential equation (ODE) $dX_t=AX_tdt, X(0)=x_0.$  Then we have $\|X(t)\| \le \|x_0\| e^{\mu(A) t }$ where 
$\mu(A)$ is the logarithmic norm of the matrix $A$ 
\cite{D59, strom75, gustaf06}. If $\mu(A) < 0$, then the ODE is asymptotically stable. 
Also, $\mu(A)$ using the matrix $2$-norm gives an estimate \cite{raha} for  the pseudospectrum  \cite{tref} of $A$:
$\max \Re \lambda_{\epsilon} (A) - \epsilon \le \mu(A)$ where $1 \gg \epsilon > 0$.
Since the pseudospectrum captures the stability of the numerical solution of the ODE over a finite number of time steps under the effect of local stiffness and nonnormality of $A$ (we shall refer to this as the numerical stability) \cite{higham, higham2}, having $\mu(A) < 0$ implies numerical stability in addition to the asymptotic stability of the ODE. As is
already shown in \cite{higham}, transient numerical stability affects the computation and choice of methods for numerical integration of the ODE.
\par
In this paper we extend the classical logarithmic norm to the stability analysis of It\^o stochastic differential equations (SDE) and introduce the stochastic logarithmic norm for estimating the numerical stability of an SDE in order to facilitate the selection of stiff and balanced stochastic numerical integration schemes. 
Letting $A, B \in \C^{n\times n}, X(t)\in \C^n,~B \in \C^{n \times n}$, the It\^o SDE with a single channel of multiplicative noise is considered in the form  of $dX_t=AX_t dt+ B X_t dW_t$ given the initial condition that $X(0)=x_0$ with probability (w.p.) $1$, where $W$ is the one dimensional 
Wiener process such that $\int_0^s dW_t \sim N(0,s)$, i.e., is standard Gaussian distributed with mean $0$ and variance $s$. 
In our definition of the stochastic logarithmic norm we shall use the the matrix $p$-norm induced by the vector $p$-norm and the expectation of the 
$l$-th raw moment, i.e., the $l$-th mean of the solution to the linear multiplicative SDE. 
The stochastic logarithmic norm is computed over the sample paths as the expected logarithmic norm of the system in the sense of the existence of a generalized derivative of the Wiener process which itself is not obligatory differentiable with respect to time. 
\par
Throughout the paper the following standard assumptions are made as in \cite{kloeden}.
Let there be a common probability space $( \Omega, {\mathcal A}, P) $ with index  
$ t \in {\mathcal T} \subset {\mathbb R} $ on which the stochastic process $X(t)$ is a collection
of random variables. The Wiener process $W = \{ W_t , t \ge t_0 \}$
is associated with an increasing family of $\sigma$-algebras $ \{  {\mathcal A}_t , t \ge t_0 \} $.
For the general case of multi-dimensional noise, each component of $\{W_t^{(i)} \}$ is $ {\mathcal A}_t $-measurable with
$ {\E} (W(t_0)) = 0$ w.p. $1$,$~{\E} (W(t) | {\mathcal A}_{t_0}) = 0$,
$~~{ \E} ((W^{(i)}_t- W^{(i)}_s)(W^{(j)}_t-W^{(j)}_s) | {\mathcal A}_s) = \delta_{i,j} (t-s) $ 
for $t_0 \le s \le t$  and 
$ \Delta W = W(t_{n+1}) - W(t_n) $, the component wise increments of the multi-dimensional Wiener process, are independent of each other at all points in the partition of the time interval ${\mathcal T}$: $ t_0 \le t_1 \le t_2 \ldots \le t_r \le t_{r+1} \le \ldots \le t_N = t_f $. 
The initial value $ X _{0}$ is assumed to be $\mathcal{A}_{t_0}$-measurable 
with $\left \| X _{0} \right \|_p  < \infty$ w.p. $1$. All expectations on a function $\phi(X_t)$ are 
evaluated as $\E( \phi(X_t) | {\mathcal A}_t)$ unless otherwise stated. Inequalities and equalities involving random variables hold almost surely where applicable.
\par
The stability analysis in \cite{mit02}
uses test equations with scalar and $2$-by-$2$ matrix coefficients having multiplicative
noise of dimension one \cite{mit96}. Some of the analysis uses
the classical logarithmic norm to establish the stability of the moment equations (derived from the SDE) which are deterministic.
The present approach with stochastic logarithmic norm 
generalizes the stability analysis of SDEs as found in \cite{mit96, mit02}.

\subsection{Classical Logarithmic Norm}
For $ 1 \leq p\leq \infty,$ the %H\"{o}lder's 
$p$-norm on $\C^n$ is given by 
\begin{eqnarray*}
 \|x\|_p :=
\left\{\begin{array}{ll}\left(\sum^n_{j=1}|x_j|^p\right)^{1/p}, &  \mbox{ for } 1 \leq p < \infty,\\
\max_{1\leq j\leq n}|x_j|, & \mbox{ for } p = \infty. \end{array}\right.
\end{eqnarray*}
Obviously, $\|x\|_2 := (x^H x)^{1/2}$ is the $2$-norm on $\C^n$.
For $ A \in \cnn,$ the spectrum $\Lam(A)$ of $A$ is given by $\Lam(A) :=\{ \lam \in \C : \rank(A-\lam I) < n\}.$  
We denote a matrix $p$-norm on $\C^{n\times n}$  induced by the vector $p$-norm as $\|\cdot\|_p$ for $p=1,2, \infty$ and define these norms as
$\|A\|_2 :=\max_{j}\{\sqrt{\lam_j} : \lam_j\ \in \Lambda(AA^H)\} ,$
$\; \|A\|_1:=\max_{j=1, \dots, n}\left(\sum_{i=1}^n|a_{ij}|\right)\, \mbox{and}\, 
\|A\|_{\infty}:=\max_{i=1, \dots, n}\left(\sum_{j=1}^n|a_{ij}|\right)$. Then the 
logarithmic norm for a single matrix is defined as $$\mu_p(A) := \lim_{h\rightarrow 0^{+}}\frac{\|I+ hA\|_p-1}{h}.$$ 
For the $1, 2$ and $\infty$-norms, the classical 
logarithmic norms, respectively, are computed  (\cite{hairer97},Vol 1) as
$\mu_{1}(A):=\max_{j}\left(\Re(a_{jj})+\sum_{j\not=i}^n|a_{ij}|\right),\,\mu_2(A):=\dm{\frac{\lam_{\max}(A+ A^{H})}{2}}\,$ and as 
$\mu_{\infty}(A):=\max_{i}\left(\Re(a_{ii})+\sum_{i\not=j}^m|a_{ij}|\right)$.

\subsection{Stability of the SDE}
The stability of the vector SDE with single channel multiplicative noise  is defined as follows.
\begin{definition} \cite{burrage, burrage2}
The equilibrium solution $X_t \equiv 0$ to $dX_t=AX_t dt + B X_t dW_t$ is stochastically 
stable in the $l$-th mean ($l$ is a finite integer $\ge$ 1) using a vector $p$-norm  if  $\forall \epsilon>0$, $\exists
\delta > 0 $ such that
\begin{equation}
\E(\|X(t)\|_p^l)<\epsilon \quad \forall\, t \ge t_0 \quad and  \quad \|X(t_0)\|_p < \delta ~w.p.~1
\end{equation}
and  is asymptotically stable in $l$-th mean if in addition, $\exists \delta_0 >0$ such that
\begin{equation}
\lim_{t\to\infty}\E(\|X_t\|_p^l)=0 \quad \forall\,\, \|X(t_0)\|_p < \delta_0 ~w.p.~1.
\end{equation}
\end{definition}

\section{Background} \label{prelim}
In the following items we review  the existing stability analysis of 
linear stochastic differential equations with multiplicative noise.
\begin{itemize}
\item[$(a)$] In \cite{mit96} scalar stochastic differential equations of the form 
$dX_t=\lam X_t dt+\beta X_t dW_t$ with  $X_0=1$ w.p. $1$,  where $\lam$ and $\beta$ are constants 
have been considered and it is shown that the above stochastic differential 
equation is mean square stable using the $2$-norm if $2\Re(\lam)+|\beta|^2 \leq 0$ when $\lam, \beta$ are complex scalars.
\item[$(b)$] In \cite{mit02} vector stochastic differential equations, with single channel multiplicative noise, of the form 
$dX_t= D X_t dt+ B X_t dW_t,$ where $D=\left( \begin{matrix}{\lam_1 & 0 \cr 0 & \lam_2} \end{matrix} \right)$ and 
$B=\left( \begin{matrix}{\alpha_1 & \beta_1 \cr \beta_2 & \alpha_2} \end{matrix} \right)$ have been analyzed for mean
square stability. It is shown that the above stochastic differential equation is mean square stable in the $\infty$ 
norm if $\max\{2\lam_1 + (|\alpha_1|+|\beta_1|)^2, 2\lam_2 + (|\alpha_2|+|\beta_2|)^2\}<0.$
\item[$(b_1)$] In \cite{mit96} for the same stochastic differential equation as in $(a)$ 
the Mil'stein scheme \cite{mile74} $X_{n+1} = X_n +\lam X_n h+\mu X_n \Delta W + \frac{(\Delta W)^2-h}{2}
\mu^2 X_n $ with $O(h^{1.5})$ root mean square error is analyzed for stability in the mean square sense. 
The mean square stability function $R(h)=|1+h\lam|^2+|h\mu^2|+\frac{1}{2}|h^2\
\mu^4|$ is derived and it is shown that the SDE is stochastically asymptotically stable in the mean square when $R(h)<1$.
\item[$(b_2)$] The reference \cite{mit96} also considered the SDE as given in $(b)$ and applied the Euler-Maruyama numerical integration scheme as $X_{n+1}=X_n + h DX_n+ B X_n\Delta W,$ where $h$ and $\Delta W$ stand for step-size and the increment 
of the Wiener process, respectively. Then the discretized SDE is shown to be stochastically asymptotically stable if $$\max\{(1+\lam_1 h)^2 + 
(|\alpha_1|+|\beta_1|)^2, (1+\lam_2 h)^2 + (|\alpha_2|+|\beta_2|)^2\} < 1.$$
\end{itemize}
\par
%Before defining the stochastic logarithmic norm, 
Next we briefly review a few essential aspects of the logarithmic norm that are used in the stability analysis of ODEs.
\begin{itemize}
\item[$(c)$] In \cite{D59, lozi58} the logarithmic norm was introduced in order to derive error bounds 
for the solution of initial value ODE problems using differential inequalities that distinguish between forward and reverse time integration. This led to requirements for the stability of initial value and boundary value ODE problems. 
The classical analysis, using the vector norm $\norm{}$ on $\C^n$ and the sub-ordinate matrix norm on $\C^{n\times n}$,
defined the logarithmic norm of a matrix $A$ as $$\mu_p(A):=\lim_{h\rightarrow 0^{+}}\frac{\|I+ h A\|_p-1}{h}.$$
\item[$(d)$] More recently S\"oderlind  \cite{gustaf06} considered $f: D\subset X\rightarrow X$ and defined least upper bounds (lub) and 
greatest lower bounds (glb) Lipschitz constants by
$$L[f]=\sup_{u\not = v}\frac{f(u)-f(v)}{|u-v|},~~ l[f]=\inf_{u\not= v}\frac{f(u)-f(v)}{|u-v|},$$ for $u$ and $v\in D,$ where 
the domain is path connected and $$l[f] |u-v|\leq |f(u)-f(v)|\leq L[f]|u-v|.$$  If $l[f]>0,$ then $f(u)\rightarrow f(v)$ 
implies that $u\rightarrow v.$ Then $f$ is an injection, with an inverse on $f(D),$ the same as a matrix 
$A\in \C^{n\times n}$ 
being invertible if its glb is strictly positive. Also then, $L[f^{-1}]=l[f^{-1}],$ where $L[f^{-1}]$ is defined over $f(D).$ 
If $f=A$ is a linear map then $L[A]=\|A\|.$ Hence $L[.]$ is left and right $G$-differentiable for the class of Lipschitz maps. 
This allows one to define the lub and glb logarithmic Lipschitz constant, by $$M[f]=\lim_{h\rightarrow 0^{+}}\frac{L[I+h f]-1}{h}, m[f]=
\lim_{h \rightarrow 0^{-}}\frac{L[I+hf]-1}{h}.$$ The lub logarithmic Lipschitz generalizes the classical logarithmic norm for 
every matrix $A,$ so that, $M[A]=\mu(A).$    
\end{itemize}  
In the following section we develop the stochastic logarithmic norm as an 
upper bound estimate of the rate of growth of the solution of a multiplicative noise linear SDE. The rate of growth is analyzed as a Dini derivative of the $l$-th
mean in vector $p$-norm of the solution. This approach may be seen as a special case of the modern definition of
logarithmic norm in \cite{gustaf06} and as an extension of the classical definition as in \cite{D59, lozi58, strom75}.
Other works \cite{higueras} deal with logarithmic norm of matrix pencils with one invertible matrix. However
in our current treatment we do not use matrix pencils for defining the stochastic logarithmic norm.

\section{Definition of the Stochastic Logarithmic Norm}
A linear It\^o SDE with a single channel of multiplicative noise is written as
\begin{equation}\label{lsde}
dX_t =AX_t dt +B X_t dW_t,
\end{equation}
where $A,B\in \C^{n\times n},$ are constant matrices. 
For the non-linear  It\^o SDE with a single channel multiplicative noise given by 
\begin{equation}
dX_t=f(X,t) dt+g(X,t) dW_t ,\label{gensde}
\end{equation}
the strong order $1.0$ It\^o-Taylor expansion (Mil'stein scheme when used in numerical integration) of $X_t$ at $t=t_{n+1}$ is given as 
\begin{equation}
X_{n+1} = X_n + f(X_n,t_n)h+g(X_n,t_n)\Delta W +\left(g\frac{\partial{g}}{\partial{x}}\right)_{(X_n,t_n)} 
\frac{((\Delta W)^2-h)}{2} + R, \label{milnl}
\end{equation}
where $h=t_{n+1} - t_n $ is the step-size, $W(t_{n+1}) - W(t_n) =: \Delta W \sim N(0,h)$ and $R$ are the $O(h^{1.5})$ remainder terms in the root mean square sense.
When linearized at $t=t_n$, $A := \left( \frac{\partial f}{\partial x}\right)_n$ and $B:=\left( \frac{\partial g}{\partial x} \right)_n$.
We shall define the stochastic logarithmic norm 
by studying the growth of the $l$-th mean of the solution $X$ of the linear SDE (\ref{lsde}) 
(linearized SDE in case of SDE (\ref{gensde})) w. r. t. time. 
We seek an upper bound for the  growth rate of $\E(\|X\|_p^l)$ using the upper-right Dini derivative which for any
function $\Xi(t)$ w.r.t $t$ is defined as
\begin{equation}\label{dini}
D_+\left(\Xi(t)\right)=\lim_{h\to0+}\frac{\Xi(t+h)-\Xi(t)}{h}.
\end{equation}
\par
The strong order $1.0$ It\^o-Taylor expansion (linearized at $t$ in case of the non-linear SDE (\ref{gensde})) applied to $X_t$ over the time interval  $[t,t+h]$ gives
\begin{equation}
X(t+h)=X(t)+hAX(t)+\Delta WBX(t) + \frac{(\Delta W)^2-h}{2}B^2X(t) + R, \label{millinear}
\end{equation}
with the remainder terms $R$ being $O(h^{1.5})$ in the root mean square sense.  
In the above expansion (\ref{millinear}) $\Delta W:=\left\{\Delta W(t), t\geq 0\right\}$ are the 
independent increments of a Wiener process over the interval $[t,t+h]$. 
Taking the $p$-norms and raising both sides of (\ref{millinear}) to the power of $l$, we can write the following inequality:
\begin{equation}
\|X(t+h)\|_p^l \le \| I+hA+\Delta WB+\frac{(\Delta W)^2-h}{2}B^2 + R_x \|_p^l \|X(t)\|_p^l,
\end{equation}
where $R_x$ are root mean square $O(h^{1.5})$ remainder terms.
For the $l$-th mean using a $p$-norm,  we apply expectation to both sides of the above inequality and get, almost surely (a.s.),
\begin{equation}\label{growthexp}
\E\|X(t+h)\|_p^l \le  \E\left(\|I+hA+\Delta WB+\frac{(\Delta W)^2-h}{2}B^2 +R_x \|_p^l \right) \E\|X(t)\|_p^l.
\end{equation}
observing that $X(t)$ is independent of the Wiener increment $\Delta W$ since the Wiener process is a non-anticipative process.
Then,  we estimate the expected rate of growth as
$$ \E \left( D_+ \| X(t)  \|_p^l \right)  \le  \lim_{h \to 0^+} \frac{ \E( \| I +hA+\Delta WB+\frac{(\Delta W)^2-h}{2}B^2 + R_x \|_p^l )  - 1 }{h} \E(\|X(t)\|_p^l ).$$
Based on the limit term on the right hand side above we
 introduce the stochastic logarithmic norm.
\begin{definition} The {\bf stochastic logarithmic norm} of a square matrix pair of same dimensions, $(A,B)$, in the $l$-th mean using a matrix $p$-norm is defined as
\begin{equation}\label{nudef}
\nu_p^l\left(A,B\right)=\lim_{h\to0^{+}}\frac{\E\left(\|I+hA+\Delta WB+\frac{(\Delta W)^2-h}{2}B^2 \|_p^l\right)-1}{h},
\end{equation}
where the limit is taken in the sense of the existence of the generalized derivative of the Wiener process and 
$\|A\|_p,\|B\|_p$ are assumed to be finite.
\end{definition}
From the above definition the expected rate of growth of the solution can be estimated as $ \E \left( D_+ \| X(t)  \|_p^l \right)  \le \nu_p^l\left(A,B\right) \E(\|X(t)\|_p^l )$ since $R_x$ are of root mean square $O(h^{1.5})$. Obviously the weaker estimate $$ D_+\E \| X(t)  \|_p^l   \le \nu_p^l\left(A,B\right) \E(\|X(t)\|_p^l )$$ also
holds so that $\E \| X(t)  \|_p^l  \le e^{ \nu_p^l\left(A,B\right) t} \|X(t_0)\|_p^l$ for the linear SDE (\ref{lsde}). 
The linear SDE (\ref{lsde}) is stochastically stable in the $l$-th mean using a $p$-norm when $\nu_p^l\left(A,B\right) \le  0$
and is asymptotically stable in the $l$-th mean using a $p$-norm when $\nu_p^l\left(A,B\right) < 0$.
The linearized SDE (\ref{gensde}) (and hence the linear SDE (\ref{lsde})) is numerically stochastically stable in the $l$-th mean using the $p$-norm if $\nu_p^l \le 0$ since over $k ~(k \ll \infty)$  time steps each of size sufficiently small $h_i>0$ almost surely we have $\E\|X_{t+kh}\|_p^l \le e^{\sum_{i=1}^k \nu_p^l(A_i, B_i) h} \E\|X_{t}\|_p^l$ when each $\nu_p^l(A_i, B_i) \le 0$. If $\E\|X_{t}\|_p^l$ is finitely bounded w.p. $1$, then  $\E\|X_{t+kh}\|_p^l$ is almost surely finitely bounded. 
Later in the paper we relate the stochastic logarithmic norm to the pseudospectrum of the Ito stability matrix $A-\frac{1}{2}B^2$ for the linear SDEs with multiplicative noise.
\par
It may be remarked that for an SDE with additive noise and for an ODE, the stochastic logarithmic norm is given as
$\nu_p^1(A,0) = \mu_p (A)$.
\section{Some Estimates of the Stochastic Logarithmic Norm}\
In this section we derive some estimates of the stochastic logarithmic norm using the matrix  
$p$-norm, where $p=1, 2, \infty.$ The estimates show the incremental behavior of the stochastic logarithmic norm
under perturbations and also the effect of noise on the deterministic ODE. \par
Let $\lambda_{\max}(A)$ be 
the largest eigenvalue of any square matrix $A$.
We state the following Lemma from pp.62, \cite{bhatia97} which we shall use in estimating the stochastic 
logarithmic norm while using the matrix $2$-norm.
 \begin{lemma} \label{lemamain}
%\begin{itemize}
%\item[$(a)$] 
Let $A, B\in \C^{n\times n}$ be Hermitian matrices. Then
$$\lam_{\max}{(A)} + \lam_{\min}{(B)} \leq \lam_{\max}(A + B)\leq \lam_{\max}(A)+ \lam_{\max}(B).$$
The equality holds when $B=k I,$ where $k$ is a scalar constant and $I$ is an identity matrix. 
\end{lemma}
For the matrix $2$-norm and in the $l^{th}$ mean we can estimate the following from (\ref{nudef}).
\begin{theorem}\label{main12}
 The stochastic logarithmic norm in the $l^{th}$ mean, where $1 \leq l < \infty$ using the matrix $2$-norm satisfies the 
following bounds in the sense of existence of a generalized derivative for the Wiener process as $\zeta dt = d W_t,~ \zeta \sim N(0,1)$.
\begin{eqnarray}
\nu_2^l(A,B)  & \le & \frac{l}{2}\lambda_{\max}(A+A^H)+ \frac{l}{4} \left( \lambda_{\max}(B+B^H) + \lambda_{\max}(-B-B^H) \right) \nonumber \\ && + \frac{l}{2}\lambda_{\max}(B^H B)+
\frac{l(l-2)}{8}\lambda_{\max}^2(B+B^H), ~~ l > 2, \label{lambdanu2est} \\
\nu_2^l(A,B) & \le & \frac{l}{2}\lambda_{\max}(A+A^H)+ \frac{l}{4} \left( \lambda_{\max}(B+B^H) + \lambda_{\max}(-B-B^H) \right) \nonumber \\ && +\frac{l}{2}\lambda_{\max}(B^H B),~l \le 2, \label{lamnu2est12} \\
\nu_2^l(A,I) & = & \frac{l}{2}\lambda_{\max}(A+A^H)+\frac{l}{2}+\frac{l(l-2)}{2}, \label{lambdanu21est} 
% \nu_2^1(A,I) & = & \frac{1}{2}\lambda_{\max}(A+A^H)  = \mu_2(A) \label{lambdanu21est12} \\
% \nu_2^2(A,I) & = & 2 \mu_2(A)  + 1 \label{lambdanu22est} 
\end{eqnarray}
where $A,B \in {\mathbb R}^{n \times n}$ are the (linearized) drift and diffusion coefficient matrices in a vector SDE with
a single channel of Wiener process noise.
\end{theorem}
\pf For $1\leq l <\infty$ we have from (\ref{nudef})
$$\nu_2^l\left(A,B\right)=\lim_{h\to 0^{+}}\frac{\E\left(\|I+hA +\Delta W B+\frac{(\Delta W)^2-h}{2}B^2\|_2^l\right)-1}{h}=$$
$$\lim_{h\to 0^{+}}\frac{\E\left(\sqrt{\lambda_{\max}\left((I+hA+\Delta W B+\frac{(\Delta W)^2-h}{2}B^2)^H (I+h A+\Delta W B+
\frac{(\Delta W)^2-h}{2}B^2)\right)}\right)^l}{h}.$$
Then we may write
\begin{eqnarray*}
&&(I+hA+\Delta W B+\frac{(\Delta W)^2-h}{2}B^2)^H(I+h A+\Delta W B+
\frac{(\Delta W)^2-h}{2}B^2)\\&&= (I+h A^H +\Delta W B^H + \frac{(\Delta W)^2-h}{2}(B^H)^2)(I+hA + \Delta W B + 
\frac{(\Delta W)^2-h}{2}B^2)=\\
&& I + h A + \Delta W B + \frac{(\Delta W)^2-h}{2}B^2 + h A^H + \Delta W B^H + 
(\Delta W)^2 B^H B \\ &&+\frac{((\Delta W)^2-h)(B^H)^2}{2}+ \cdots
= I + h(A+ A^H) + \Delta W (B+B^H)+(\Delta W)^2 B^H B+\\ && \frac{(\Delta W)^2-h}{2}(B^2 + (B^H)^2)+ \cdots.
\end{eqnarray*}
It is possible to write
$C= h(A+ A^H) + \Delta W (B+B^H)+(\Delta W)^2 B^H B+ \frac{(\Delta W)^2-h}{2}(B^2 + (B^H)^2)$ in a series in
the normalized time step size $h \ll 1$.
We use the identities $\lam_{\max}(c A)= |c| \lam_{\max}(\sign(c) A)$ (where $c$ is a constant, $\sign(c)=\dm{\frac{z}{|z|}}, z\not=0$) and  
$\lambda_{\max}^2(A)=\lambda_{\max}(A^2)$ for estimating the terms in the series.
Using the triangle inequality
and the identity $[\lam_{\max}(I + C)]^{\frac{l}{2}}= (1 + \lam_{\max}(C))^{\frac{l}{2}}$,  the following can be written:

\begin{eqnarray*}
&&(1+ \lam_{\max}(C))^{\frac{l}{2}} = 1+\frac {l}{2}\lam_{\max}(C)+\frac{l(l-2)}{8}\lam_{\max}(C)^2+\dots = 1+\\
&& \frac{l}{2}\lam_{\max}\left( h(A+ A^H) + \Delta W (B+B^H)+(\Delta W)^2 B^H B+ \frac{(\Delta W)^2-h}{2}(B^2 + (B^H)^2)\right)+ \\
\end{eqnarray*} 
$$\frac{l(l-2)}{8}\lam_{\max}\left(h(A+ A^H) + \Delta W (B+B^H)+(\Delta W)^2 B^H B+ \frac{(\Delta W)^2-h}{2}(B^2 + (B^H)^2)\right)^2$$
\begin{eqnarray*}
&& +\cdots \label{2estneq}\leq 1+ \frac{l}{2}h\lam_{\max}(A + A^H) +\frac{l}{2}  \lam_{\max}(\Delta W (B + B^H)) +\frac{l}{2}(\Delta W)^2
\lam_{\max}(B^H B)
\end{eqnarray*}
\begin{eqnarray}
&& +((\Delta W)^2-h)\frac{l}{4} \lam_{\max}(B^2 + (B^H)^2)+(\Delta W)^2\frac{l(l-2)}{8}\lam_{\max}^2(B+ B^H)+\cdots.\label{2esteq}
\end{eqnarray}
The equality holds when $B=I$ and follows from Lemma \ref{lemamain}.  
Taking expectation on both sides and  in the sense of generalized derivative of the Wiener process which is the Gaussian white noise $\zeta$ 
 we have 
$ \E(\lam_{\max}(I + C))^{\frac{l}{2}}=\E(1+ \lam_{\max}(C))^{\frac{l}{2}} \leq 1+ \frac{l}{2}h\lam_{\max}(A + A^H) + \E(|\zeta|) h
\frac{l}{4}( \lam_{\max}(B + B^H) + \lam_{\max}(-B -B^H) )+\E((\Delta W)^2)\frac{l}{2} \lam_{\max}(B^H B)
+(\E((\Delta W)^2)-h) \frac{l}{4}\lam_{\max}(B^2 + (B^H)^2)+\E((\Delta W)^2)\frac{l(l-2)}{8}\lam_{\max}^2(B+ B^H) + O(h^{1.5})
= 1+h \frac{l}{2}\lam_{\max}(A + A^H)+ h \frac{l}{4}( \lam_{\max}(B + B^H) + \lam_{\max}(-B -B^H) )  + \frac{l}{2}h\lam_{\max}(B^H B)-h\frac{l(l-2)}{8}\lam_{\max}^2(B + B^H)+O(h^{1.5}).$ 
For $h\to 0^{+}$ we obtain
$$\lim_{h \rightarrow 0^{+}}\frac{\E(\lam_{\max}(1+ C))^{\frac{l}{2}}-1}{h} \leq \frac{l}{2}\lam_{\max}(A + A^H)+
\frac{l}{4}( \lam_{\max}(B + B^H) + \lam_{\max}(-B -B^H) )$$ $$ + 
\frac{l}{2}\lam_{\max}(B^H B)+\frac{l(l-2)}{8}\lam_{\max}^2(B + B^H).$$ 
{\vskip 1 ex} \noindent $(b)$ When $ B=I,$ we have 
\begin{eqnarray*}
\lim_{h \rightarrow 0^{+}}\frac{\E(\lam_{\max}(I+ C))^{\frac{l}{2}}-1}{h} &=& \frac{l}{2}\lam_{\max}(A + A^H)+\frac{l}{2}\lam_{\max}(B^H B)+
\\ && \frac{l(l-2)}{8}\lam_{\max}^2(B + B^H)\\ \Rightarrow  \nu_2^l(A,B) &=& \frac{l}{2}\lambda_{\max}(A+A^H)+\frac{l}{2}+\frac{l(l-2)}{2}. \blacksquare 
\end{eqnarray*}
\begin{corollary}
For a positive integer $l > 2$, 
\begin{eqnarray}
\nu_2^l(A,B) & \le & l \left ( \mu_2 (A) + \frac{1}{2} \| B \|_2^2 + \frac{1}{2} (\mu_2(B)+\mu_2(-B)) + \frac{l-2}{2} ( \mu_2(B) )^2 \right ) \label{msest}  
%\nu_2^l(A,B) & \le & l \mu_2(A - \frac{1}{2} B^2) +  \frac{1}{2} \| B \|_2^2 + \frac{l}{2} (\mu_2(B)+ \mu_2(-B)) \nonumber \\
%&& + \frac{l}{2} \| B \|_2^2 + \frac{l}{2} \mu_2(B^2) + \frac{l(l-2)}{2} (\mu_2 (B))^2 \label{p2muest}
\end{eqnarray}
For $l \le 2 $, 
\begin{eqnarray}
\nu_2^l(A,B) & \le & l \left( \mu_2 (A) + \frac{1}{2} (\mu_2(B)+\mu_2(-B)) + \frac{1}{2} \| B \|_2^2 \right)  \label{ms12est}  
\end{eqnarray}
\end{corollary}
\pf 
The results (\ref{msest}) and (\ref{ms12est}) follow directly from  (\ref{2esteq}) in Theorem \ref{main12} by
Lemma \ref{lemamain} and the It\^o isometry for the expectation of the It\^o integrals.$\blacksquare$ \\
The inequality (\ref{msest}) has been used for mean square stability estimates in \cite{mit02}.  
\begin{corollary}
If $l=1, 2$ and $p=2$, that is, in the mean and in the mean square and using the matrix two norm, the following results hold
when $B$ is the identity matrix $I$.
\begin{eqnarray}
&&  \nu_2^1(A,I)=\lambda_{\max}\left(\frac{A+A^H}{2}\right)+ \frac{1}{2}-\frac{4}{8}=\lambda_{\max}\left(\frac{A+A^H}{2}\right)=\mu_2(A), \label{beq1mean} \\ 
&& \nu_2^2(A,I)=\lambda_{\max}\left(A+A^H\right)+1 = 2 \mu_2(A)+1, \label{beq1meansq}
\end{eqnarray}
\end{corollary}
\pf  The results follow from Theorem \ref{main12}, by substituting $p=2,l=1$ and $p=2, l= 2.\blacksquare$
\par
Some properties of the stochastic logarithmic norm under perturbation is obtained in the following result.
In the generalized sense of the Wiener process derivative a couple of useful inequalities showing how
the Wiener process perturbs the deterministic ODE is also given.
\begin{theorem} For any square matrices of same dimensions, $B, A, \Delta A, \Delta B$ and a real number $\alpha > 0$  
the stochastic logarithmic norm has the following properties.
\begin{eqnarray}
\nu_p^l(\alpha A, \sqrt{\alpha} B) & =&  \alpha \nu_p^l(A, B)  \label{scal} \\
\nu_p^{1} (A+\Delta A, B+\Delta B)  & \le &  \nu_p^1(A, \sqrt{2} B) +  
           \nu_p^1(\Delta A, \sqrt{2} \Delta B) + \frac{1}{\sqrt{2}} \| (B-\Delta B)^2 \|_p  \label{tri1} \\
\nu_p^l (A+\Delta A, B+\Delta B)  & \le &  \nu_p^l (A, \frac{B+\Delta B}{\sqrt{2}} ) +  \nu_p^l (\Delta A, \frac{B+\Delta B}{\sqrt{2}} )  \label{tri2} \\
\nu_p^l(A,0)  & = & l \mu_p(A) \label{detmunu} 
\end{eqnarray}
Further, in the sense of the existence of a generalized derivative of the Wiener process which is a Gaussian white noise, the following estimate holds for a positive integer $l$ and for any of the matrix $p$-norms.
\begin{eqnarray}
\nu^l_p (A,B) & \le & l \mu_p(A) + \frac{l}{2} \mu_p(-B^2) +  \frac{l(l+1)}{4} \| B \|^2_p  + l \| B \|_p  \label{lpest1} \\
& \le & l \mu_p(A) +  l \| B \|_p \left( 1 + \frac{l+3}{4} \|B\|_p \right)  \le l \mu_p(A) + l  \left( 1 + \frac{l+3}{4} \|B\|_p \right)^2 \label{lpest} 
\end{eqnarray}
\label{propthm}
\end{theorem}
\pf We can scale the Wiener process as $\sqrt{\alpha}$ and write 
\begin{eqnarray}
\nu_p^l(\alpha A, \sqrt{\alpha} B) & =  & \lim_{h \to 0^+}  \alpha \frac{  \E \left(\| I + (\alpha h) A + \Delta (\sqrt{\alpha}W ) B +   \frac{(\Delta 
\sqrt{\alpha}W)^2 - \alpha h}{2}  B^2 \|_p^l\right)-1}{ \alpha h} \nonumber \\
& = & \alpha \nu_p^l(A,  B)  \nonumber
\end{eqnarray}
for any positive integer $l$ and any of the matrix $p$-norms.  
\par
From the definition of the stochastic logarithmic norm and scaling the Wiener process, we can write 

$$ \nu_p^1 (A+\Delta A, B+\Delta B)  \le   \lim_{h \to 0^+}\bigg(  
 \frac{  \E (  \| I + (2h) A + ( \sqrt{2} \Delta W ) ( \sqrt{2} B) + \frac{ (\Delta \sqrt{2}W)^2 -2 h}{2} ( \sqrt{2} B)^2 \|_p )-1}{2h}$$ 
$$ + \frac{\E ( \| I + (2 h) \Delta A +(\sqrt{2} \Delta W) ( \sqrt{2} \Delta B ) + \frac{( \Delta  \sqrt{2}W)^2 - 2 h}{2}  ( \sqrt{2} \Delta B)^2 \|_p )-1}{2h}+ $$
$$ \frac{\E | \int_0^h W_u dW_u |}{h} \|B \Delta B + \Delta B B - B^2 - (\Delta B)^2 \|_p  \bigg) \\
 \le   \nu_p^1 (A,\sqrt{2} B) + \nu_p^1(\Delta A, \sqrt{2} \Delta B)  + \frac{1}{\sqrt{2}} \| (B-\Delta B)^2 \|_p  
$$
since $\lim_{h \to 0^+} \frac{\E | \int_0^h W_u dW_u |}{h}= h/(\sqrt{2}h) = \frac{1}{\sqrt{2}}$.  For any positive integer $l$, the above can be re-written as

\begin{eqnarray*} 
\nu_p^l (A+\Delta A, B+\Delta B) & \le &  \lim_{h \to 0^+} \bigg(  
  \bigg( 2^{l-1} \frac{1}{2^{l-1}} \E \| I + (2h) A + ( \sqrt{2} \Delta W ) ( \sqrt{2} \frac{B+\Delta B}{2}) \\ 
&& + \frac{ (\Delta \sqrt{2}W)^2 -2 h}{2}( \sqrt{2} \frac{B+\Delta B}{2})^2  \|_p^l-1 \bigg) / (2h)  \\
&& + \bigg( 2^{l-1} \frac{1}{2^{l-1}} \E \| I + (2h) \Delta A + ( \sqrt{2} \Delta W )( \sqrt{2} \frac{B+\Delta B}{2}) \\ 
&& + \frac{( \Delta \sqrt{2}W)^2 -2 h}{2} ( \sqrt{2} \frac{B+\Delta B}{2})^2  \|_p^l -1 \bigg) / (2h) \bigg)  \\
&& \le   \nu_p^l (A,\frac{B+\Delta B}{\sqrt{2}}) + \nu_p^l(\Delta A, \frac{B + \Delta B}{\sqrt{2}}). 
\end{eqnarray*}
\par
Applying (\ref{tri1}) recursively to $\nu_p^l(A+0,B/2 + B/2)$ and using (\ref{scal}), we obtain $\nu_p^l(A,B) \le \lim_{n \to \infty} 
( \nu_p^l(A, \frac{1}{(\sqrt{2})^n} B ) + (2^n -1) \nu_p^l(0, \frac{1}{(\sqrt{2})^n} B ) ) =\nu_p^l(A, 0) +  \lim_{n \to \infty}(2^n/2^n) \nu^l_p(0, B) + 0 = \nu_p^l(A,0) + \nu_p^l(0,B)$. For $l=1$, the inequality reduces to  $\nu_p^1(A,B) \le \mu_p(A) + \nu_p^1(0,B)$.
\par
For the deterministic case with no noise, we have 
\begin{eqnarray*}
 \nu_p^l(A,0) & = & \lim_{h \to 0^+} \frac{\| I + hA \|^l_p - 1}{h} \\
& = & \lim_{h \to 0^+} l \| I + hA \|^{l-1}_p  D_{+,h} \| I + hA \|_p \\
& = & \lim_{\epsilon \to 0^+} l \frac{\| I + \epsilon A \|_p - 1 }{\epsilon} \\
& = & l \mu_p(A).
\end{eqnarray*}
For the stochastic logarithmic norm of the diffusion coefficient $B$ it is possible to write 
\begin{eqnarray*}
 \nu_p^l ( 0,B )  & = &  \lim_{h \to 0^+}  \frac{ \E \| I + \Delta W B + \frac{(\Delta W)^2 - h}{2} B^2 \|_p^l - 1}{h}  \\
& \le & \lim_{h \to 0^+} \left( \frac { \| I - (2h) \frac{1}{2} B^2 \|_p^l -1 }{(2h)} + \frac{ \E \| I + 2 B \Delta W  + B^2 (\Delta W)^2 \|^l_p -1 }{(2h)}  \right) \\
& \le & \nu_p^l ( -\frac{1}{2} B^2, 0 ) +  \frac{l}{2} \| B \|^2_p + l \| B \|_p + \frac{l(l-1)}{4} \| B \|^2_p 
\end{eqnarray*} 
in the sense of the existence of a generalized derivative of the Wiener process so that $ \int_0^h \xi dt =\int_0^h dW_s $ where $\xi \sim N(0,1)$ 
is a Gaussian white noise. Further we have $\nu_p^l (-\frac{1}{2}B^2, 0 )\\ =  \frac{l}{2}\mu_p (-B^2)$ and $ \frac{l}{2}| \mu_p (-B^2) | 
\le  \frac{l}{2} \| B \|_p^2$ from the properties of the logarithmic norm \cite{hairer97}. 
Hence the estimates (\ref{lpest1}) and (\ref{lpest}) .
$\blacksquare$  
\par We remark that the stochastic logarithmic norm does not satisfy the triangle inequality 
property nor the that of the multiplication by a scalar in the same way as the (deterministic) logarithmic norm. 
However, it is consistent with It\^o calculus
in its property (\ref{scal}) of multiplication with a scalar.  Again, It\^o calculus makes the noise ``redistribute`` for
any additive perturbation to the drift and diffusion coefficients as found in (\ref{tri2}). \par
The logarithmic norm coincides with the stochastic logarithmic norm in special cases.
One such case is obtained from Theorem \ref{main12} by setting $l=1$ and $B=I$ (and also by setting $B=0$ in the trivial case) for $p=2$. 
\begin{theorem}\label{lem3}
$\max_{\lam} \Re(\lam(A)) \leq  \frac{1}{2} \nu_2^2(A, I) - \frac{1}{2}$ where $I$ is the identity matrix.
\end{theorem}
\pf  From (\ref{beq1meansq}) we have
$\nu_2^2(A, I)= 2 \mu_2(A) + 1 $. The logarithmic norm has the lower bound property \cite{hairer97}: $\max_{\lam} \Re(\lam(A)) \leq \mu_2(A)$. Hence the inequality. $\blacksquare$
\section{Conditions for  Mean and Mean Square Stability}
We recall that the linear SDE (\ref{lsde}) is stochastically stable in the $l$-th mean using a $p$-norm if 
$\nu_p^l(A, B) \le 0.$ To estimate of the stability of a linear SDE one of the upper bounds derived in
the last section may be used, especially,
when it is needed to compute incrementally the effect of adding noise to an ODE or an SDE with known stability estimates.  
In practice, when using the upper bounds for estimating the stability of an SDE incrementally, a small
positive number (determined by the stochastic stability region of the numerical integrator) is used as cut-off rather than a very small tolerance or zero so that the effect of non-normality and
stiffness \cite{higham} in the SDE's (both linear and locally linearized) 
transient numerical behavior in the stochastic logarithmic norm 
is taken into account.  
\par
The following results are a couple of special cases of the stochastic logarithmic norm approach to the stability of SDEs with multiplicative noise.
\begin{itemize}
\item We consider the scalar case $A=\alpha \in \C$ and $B=\beta \in \C$ in (\ref{lsde}). Then we apply Theorem \ref{main12}$(a)$ and consider
(\ref{msest}) for $l=1, 2$ the condition $\nu_{p}^l(A, B)<0,$ so that for $p=2$,  we get the condition for stochastic stability in the mean as 
\begin{eqnarray}
%\Re( \alpha ) + \frac{1}{2}(\Im(\beta))^2 \le 
\Re( \alpha ) + \frac{1}{2} |\beta|^2 \le 0 \label{meanstab2}
 \end{eqnarray}
and in the mean square as 
\begin{eqnarray} 2\Re(\alpha)+|\beta|^2 \le 0. \label{msqstab2} \end{eqnarray}  The above conditions are
essentially the same but in practice the $0$ on the right hand side is replaced by $TOL$ which is 
a small positive real number  \cite{raha} and thus $l$ becomes significant for estimates in higher moments. 
These results are the same as in review $(a)$ of Section \ref{prelim}.
\item For the matrices $A= \left( \begin{matrix}{\lam_1 & 0 \cr
     0 & \lam_2}\end{matrix} \right) $ and $B=\left( \begin{matrix}{\alpha_1 & \beta_1\cr \beta_2 & \alpha_2}\end{matrix} \right),$ the SDE (\ref{lsde})
     is mean square stable in the $\infty$-norm if 
\begin{eqnarray} 
   2 \max \{  \lambda_1  , \lambda_2 \} +  
    2 \left( \frac{5}{4} \max\{ | \alpha_1 | + | \beta_1 | , |\alpha_2 | + | \beta_2 | \}  + 1 \right)^2 \le 0. \label{stabmsqinf}
 \end{eqnarray} 
 The above condition is obtained using (\ref{lpest}) in
which 
 \begin{eqnarray*}
\nu^2_{\infty} (A,0)  & = & \lim_{h \to 0^+} \frac{ \left( \max\{ |1 + h \lambda_1|, |1 + h \lambda_2| \} \right)^2 -1 }{h} \\ 
 & = & \lim_{h \to 0^+} 2 \max\{ |1 + h \lambda_1|, |1 + h \lambda_2| \} \max\{\lambda_1,\lambda_2\} \\
 & = & 2 \max_{i=1,2} {\lambda_i} 
\end{eqnarray*} 
and $\|B\|_{\infty} =  \max\{ | \alpha_1 | + | \beta_1 | , |\alpha_2 | + | \beta_2 | \}.$
     \end{itemize}
\section{Direct Computation of the Stochastic Logarithmic Norm and Sharper Bounds}
\begin{theorem}
In the sense of the existence of a generalized derivative for the Wiener process, i.e., $ \zeta dt = d W_t, ~ \zeta \sim N(0,1)$,
we can compute the stochastic logarithmic norm directly as follows.
 $$ \nu_p^l (A,B)  = l \E \left( \mu_p \left(A - \frac{1}{2} B^2 + B \zeta \right) \right) $$
\label{compnu}
\end{theorem}
\pf 
\begin{eqnarray}
 \nu_p^l (A,B) & = & \lim_{h \to 0^+} \frac{ \E \| I + h A - \frac{1}{2} B^2 h + B  \Delta W + \frac{1}{2} B^2 (\Delta W)^2 \|_p^l - 1}{h} \nonumber \\
%& \approx &   \lim_{h \to 0^+} \frac{\E \| I + h A - \frac{1}{2} B^2 + B h \zeta  \|_p^l - 1}{h} \nonumber \\
& = &  \lim_{h \to 0^+} \E \bigg( l  \| I + h A - \frac{1}{2} B^2  h + B \Delta W + \frac{1}{2} B^2 (\Delta W)^2  \|_p^{l-1} \nonumber \\
  && \times D_{+,h} \| I + h A - \frac{1}{2} B^2 h + B \Delta W +  \frac{1}{2} B^2 (\Delta W)^2 \|_p \bigg) \nonumber \\
& = & l \lim_{\epsilon \to 0^+} \frac{\E \| I + (A - \frac{1}{2} B^2) \epsilon + B \int_0^{\epsilon} \zeta ds  + \frac{1}{2} B^2 \left( \int_0^{\epsilon} \zeta ds \right)^2 \| -1}{\epsilon} \nonumber \\
& = & l \E \left( \mu_p \left(A - \frac{1}{2} B^2 + B \zeta \right) \right), \label{nucomput1ch}
\end{eqnarray} 
where the last equality follows from considering the It\^o formula along with the generalized derivative of the Wiener
process. $\blacksquare$  \par
The above result shows that the stochastic logarithmic norm is the expected logarithmic norm behavior of the 
SDE. Consequently we can derive the following inequalities.
\begin{corollary}
\begin{eqnarray}  \nu_p^l (A,B)  & \le & l \mu_p (A) + \frac{l}{2} \left( \mu_p (-B^2) + \mu_p(B) + \mu_p(-B) \right) \label{mubnds}  \\
 \nu_p^l (A,B) &  \ge & l \mu_p(A) - \frac{l}{2}\left( \mu_p(B^2) +   \mu_p(B) + \mu_p(-B) \right) \label{mulbnd}
\end{eqnarray}
\label{ulnucor}
\end{corollary}
\pf From (\ref{nucomput1ch}) we have $ \nu_p^l (A,B)  \le  l \mu_p (A) + \frac{l}{2} \left( \mu_p (-B^2) + \E (\mu_p(B \zeta)) \right)$
in view of the triangular inequality of the logarithmic norm. Since $\zeta \sim N(0,1)$, we have $ \E(\mu_p(B \zeta)) = \E | \zeta |  (\mu_p(B) + \mu_p(-B) )/2 = \left( \mu_p(B) + \mu_p(-B) \right) / 2$. Again, $\mu_p(A) = \mu_p(A - \frac{1}{2} B^2 + B \zeta
+ \frac{1}{2}B^2 - B \zeta) \le \mu_p(A - \frac{1}{2} B^2 + B \zeta) + \mu_p( \frac{1}{2}B^2 )  + \mu_p( - B \zeta)$ so 
that $l \mu_p(A) -  l \mu_p( \frac{1}{2}B^2 ) -  l \E \mu_p( - B \zeta) \le \nu_p^l(A,B)$.  Hence the inequalities. $\blacksquare$
\begin{corollary}
 The stochastic logarithmic norm can be bounded as follows.
\begin{equation}
 \left| \nu_p^l(A,B) \right| 
\le l \left | \mu_p(A - \frac{1}{2} B^2) \right | + l \E \left | \mu_p ( B \zeta ) \right | \le l \left\| A - \frac{1}{2} B^2 \right \|_p  +l  \left \| B \right \|_p \label{boundednu}
\end{equation}
\end{corollary}
\pf By Theorem \ref{compnu} we have  $\left | \nu_p^l (A,B) \right | = l \left | \E \left( \mu_p \left(A - \frac{1}{2} B^2 + B \zeta \right) \right) \right| $.
Then, by Jensen inequality one obtains, a.s., $$ l \left | \E \left( \mu_p \left(A - \frac{1}{2} B^2 + B \zeta \right) \right) \right |  \le 
 l \E  \left | \left( \mu_p \left(A - \frac{1}{2} B^2 + B \zeta \right) \right) \right |. $$ From the triangular inequality we have  $$\left | \mu_p \left(A - \frac{1}{2} B^2 + B \zeta \right)  \right| \le  
 \left | \mu_p \left(A - \frac{1}{2} B^2 \right) \right | + \left  | \mu_p( B \zeta ) \right| $$ and the bound 
property of the logarithmic norm leads to
%\le  
%\left \| A - \frac{1}{2} B^2  + B \zeta \right \|_p$.  Then, by triangular inequality,  
$$ \left| \nu_p^l(A,B) \right| \le l \left\| A - \frac{1}{2} B^2 \right \|_p +l \E \left \| B \zeta \right \|_p  \le l \left\| A - \frac{1}{2} B^2 \right \|_p +l \left \| B \right \|_p.  \blacksquare$$
\section{Examples}
\subsection{Stabilization of an inverted pendulum}
It is well known \cite{sanzserna, enquist} 
that a vertical inverted pendulum can be stabilized around a mean vertical position by application of 
a suitable highly oscillatory excitation in the form of an appropriate noise. We compute the stochastic logarithmic norm of such a system to show how an appropriate noise stabilizes the system.
The equation of the inverted pendulum may be written as 
\begin{eqnarray}
d \theta  &= & v ~dt  + \epsilon v ~dW_t ,~ 1 \gg \epsilon > 0 \\
d v  & = &  \frac{g}{l} \theta ~dt + b \theta~ dW_t ,
\end{eqnarray}
where $g$ is acceleration due to gravity, $l>0$ is the effective length of the pendulum and $\theta$ is the small angular displacement from 
the mean verticalposition, i.e, $\theta=0$,
so that  $A= \begin{pmatrix} {0 & 1  \cr \frac{g}{l} & 0} \end{pmatrix} $ and $B =  \begin{pmatrix} {0 & \epsilon \cr b & 0} \end{pmatrix}$. It may be noted that $\mu_2(A) = \frac{1}{2} + \frac{g}{2l}$ and $\Re \lambda_{\max}(A) = \frac{g}{l} > 0$. Thus the system without the Wiener process excitation is unstable. The stochastic logarithmic norm of the Wiener process excited system may be computed as 
$$ \nu_2^2 (A,B) =  \E \left( \max \left \{ 1 + \frac{g}{l} + (b+ \epsilon) \zeta - b \epsilon , - \left( 1 + \frac{g}{l} + (b+ \epsilon) \zeta + b \epsilon \right) \right \} \right),$$ where $\zeta \sim N(0,1)$.
For stabilizing the pendulum around the mean position $\theta=0$ we require $\nu_2^2 (A,B) \le 0$ and get the condition that
\begin{equation}
b \ge \frac{1}{\epsilon} \left( 1 + \frac{g}{l} \right) 
\end{equation}
in which $\epsilon/(2 (1 + g/l) )$ can be interpreted physically as the amplitude of a  very wide band vertical excitation (as an approximation
to Wiener process) at the base of the pendulum. Since a Karhunen-Loeve expansion of the Wiener process \cite{kloeden} contains the high frequency terms, this also shows how the above result is consistent with the result (in \cite{enquist})  that
a small amplitude highly oscillatory wide band vertical excitation stabilizes a vertical inverted pendulum. 
% It may be noted that the results in this example is remarkably similar to the nonnormality issue discussed in \cite{higham2}.
\subsection{Nonnormality} From \cite{higham2} we take this linear SDE in the form of (\ref{lsde}) in which 
$A := \begin{pmatrix} {-1 & b  \cr 0 & -1} \end{pmatrix}, B:= \begin{pmatrix} {0 & \sigma  \cr -\sigma & 0} \end{pmatrix}, b \in \R$. Obviously the system without any Wiener process excitation, i.e, the deterministic ODE is asymptotically stable but
$\mu_2(A) = \max \{ \frac{b}{2} -1, -\left( \frac{b}{2}+1 \right) \}$  due to non-normality of $A$ and the ODE system tends to be numerically unstable in its transient behavior when $|b| > 2$. 
Formulating the system as a multiplicative noise SDE we compute the stochastic logarithmic norm
directly from (\ref{nucomput1ch}) in the mean square and using the $2$-norm as in \cite{higham2}:
$$ \nu_2^2 = \max\{ \sigma^2 - 2 \pm b \} $$ whence 
it is required that $\sigma^2 \le \min\{2 \mp b\}$ so that $\nu_2^2 \le 0$ for the mean square stability of  the SDE. If $\sigma \in {\mathbb R}$ and $|b| >2$, then it is not possible to 
numerically stabilize the SDE in the mean square by choosing an appropriate $\sigma$.
For $|b| > 2$ and $\sigma \in {\mathbb C}$, $ \sigma^2 \le 2 -|b|$ would be sufficient for exploiting the noise towards   
 numerically stabilizing the SDE system in the mean square. If $\sigma,b \in \R$ and $\sigma :=b^{-\frac{1}{4}}$, then the SDE system is stochastically stable when $1 \ge b  \ge \frac{3-\sqrt{5}}{2}$.
\subsection{Numerical examples}
\begin{exam}\label{ex1}
We consider
\begin{itemize}
\item[$(a)$]$A=\left( \begin{matrix}{-100 & 0\cr 0 & -200}\end{matrix} \right), 
B=\left( \begin{matrix}{5 & 0\cr 0 & 6}\end{matrix} \right).$
\item[$(b)$]$A=\left( \begin{matrix}{-100 & 0\cr 200 & -200 }\end{matrix}\right), B=\left(\begin{matrix}{5 & 2\cr 0 & 6}\end{matrix} \right).$
\item[$(c)$] $A=\left( \begin{matrix}{-100 & 20 \cr 0 & -200}\end{matrix} \right), 
B=\left( \begin{matrix}{5 & 2\cr 0 & 6}\end{matrix} \right).$
\item[$(d)$]$A=\left(\begin{matrix}{-100+20i &  0\cr 2 & -200+i}\end{matrix}\right), 
B=\left( \begin{matrix}{5+i  & 0   \cr 2i   & -6-10i}\end{matrix} \right).$
\item[$(e)$] $A=\left(\begin{matrix}{-100 & 20 \cr 7 & -200}\end{matrix}\right), 
                   B=\left(\begin{matrix}{5 & 2\cr 4 & 6}\end{matrix} \right).$
\item[$(f)$] $A=-100, B=10$ (\cite{mit96}).
\item[$(g)$] $A := \left(\begin{matrix}{ A_1 & A_{12} \cr 0 & A_2} \end{matrix} \right); B:= \left(\begin{matrix}{ B_1 & B_{12} \cr 0 & B_2} \end{matrix} \right)$ where $A_{1} := \left(\begin{matrix}{
  0.1 & 4 & 20 \cr
  0 & 0.1 & 5 \cr
  0 & 0 & 0.1 }\end{matrix}\right),\\
  A_{2} :=\left(\begin{matrix}{
  -0.2 & 3 & 100 \cr
  0 & -0.2 & 50 \cr
  0 & 0 & -0.2 }\end{matrix}\right), B_{1} :=\left(\begin{matrix}{
  2 & 30 & 10 \cr
  0 & 2 & 50 \cr
  0 & 0 & 2 }\end{matrix}\right),  B_{2} :=\left( \begin{matrix}{
  4 & 6 & 20 \cr
  0 & 4 & 40 \cr
  0 & 0 & 4 }\end{matrix}\right),$ \\
 $A_{12} = \left(\begin{matrix}{
  2.2857 \times 10^{-2} & -2.3547 \times 10^{-2} & -6.8279 \times 10^{-2} \cr
  9.3914 \times 10^{-2} & -9.6719 \times 10^{-2} & -2.8049 \times 10^{-1} \cr
  2.8585 \times 10^{-1} & -2.9443 \times 10^{-1} & -8.5382 \times 10^{-1} }\end{matrix}\right),$\\
and $B_{12} = \left(\begin{matrix}{
  1.2606 \times 10^{-1} & -4.6007 \times 10^{-1} & 7.0963 \times 10^{-3} \cr
  1.8156 \times 10^{-1} & -6.6259 \times 10^{-1} & 1.0235 \times 10^{-2} \cr
  1.4481 \times 10^{-1} & -5.2845 \times 10^{-1} & 8.1625 \times 10^{-3} }\end{matrix}\right).$
\item[$(h)$] $A=100 \times C_{100\times100}, B=100 \times D_{100 \times 100};~C_{ij},D_{ij} \sim U(0,1)$.
\item[$(i)$]$\left(\begin{matrix}{-100 & 0\cr 0 & -1}\end{matrix}\right).$
$B=\left(\begin{matrix}{0 & 2\cr 2 & 0}\end{matrix}\right)$ \cite{mit02}
% B=\left(\begin{matrix}{0  & \sigma  \cr -\sigma & 0}\end{matrix}\right), \mbox{where}\,b=2\times 10^{11}\, \mbox{and}\, \sigma=b^{-1/4}.$ \cite{higham2}
\end{itemize} 
\begin{table}[hbtc]
\begin{center}
\begin{tabular}{|c|c|c|c|c|c|c|}
  \hline
Example &         $\mbox{Lbound}$ &            $\nu_2^2(A, B)$          & $\mbox{Ubound}$ \\
  \hline
$(a)$ &  $-1.1239\times 10^{2} $ &            $-1.0470\times 10^{2}$                & $-4.0393\times 10^{1}$    \\
   \hline
$(b)$ &  $-1.1919\times 10^{2} $  &           $-1.1468\times 10^{2}$               & $-3.1393\times 10^{1}$     \\
 \hline
$(c)$ &  $ -2.4082\times 10^{2}  $  &         $-2.2415\times 10^{2}$           & $-1.5302\times 10^{2} $     \\
 \hline
$(d)$ & $-2.2490\times 10^{2} $ &              $ -2.2354\times 10^{2}$          & $-5.9075\times 10^{1} $   \\
 \hline
$(e)$  & $ -2.6837\times 10^{2}$  &               $-2.3232\times 10^{2}$           & $ -1.21915\times 10^{2}$   \\
\hline
$(f)$   & $-3.0000\times 10^{2} $ &             $ -3.0026\times 10^{2}$           & $-1.0000\times 10^2$   \\
 \hline
$(g)$  &  $-9.1852\times 10^{2} $ &             $9.2453\times 10^{2}$            &$4.8398\times 10^{3}$    \\
\hline
$(h)$  & $-2.5191\times 10^{7}$  &              $1.2369\times 10^{5}$             & $2.5330\times 10^{7}$ \\
\hline
$(i)$ & $-6.0000 $               &               $-5.91409$         & $-2.0000$             \\ 
% \hline  $(j)$ &  $-6.0000$               &               $ -5.9198 $        &          $2.0000$  \\
\hline

\end{tabular}
\caption{ The table gives values $\nu_2^2(A, B)$ and compares them with the estimates $Ubound$ and $Lbound$.}\label{T1}
\end{center}
\end{table}
In Table \ref{T1} Ubound is computed using the right hand side in (\ref{mubnds}) with $p=2,~l=2$ and 
Lbound is computed using the right hand side in (\ref{mulbnd}). The stochastic logarithmic norm is
computed using (\ref{nucomput1ch}).
\end{exam}
\section{Extension to Multiple Noise Channels}
The definition of the stochastic logarithmic norm can be extended to the vector SDE with multiple channels of multiplicative noise (i.e., with multi-dimensional Wiener process). Considering the It\^o-Taylor strong order 1.0 (Mil'stein) scheme for the SDE
\begin{eqnarray}  dX_t = A X_t dt + \sum_{j=1}^{m} B^{(j)} X_t dW_t^{(j)}, \label{mSDE} \end{eqnarray} 
where each $W^{(j)}$ is an independent component Wiener process, we can define the following.
\begin{definition}
The {\bf stochastic logarithmic norm} for a tuple of $m+1$ square matrices of same dimensions, $(A, B^{(1)}, B^{(2)}, \cdots, B^{(m)})$, i.e., $\left(A, B^{(1:m)} \right)$, which are the (linearized) drift and diffusion coefficients of a (non-linear) vector SDE  with $m$ channels of multiplicative noise, is defined as 
$$ \nu_p^l \left(A, B^{(1:m)} \right) =
\lim_{h \to 0^+} \frac{\E \| I + hA + \dm{\sum_{j=1}^m B^{(j)} \Delta W^{(j)} + \sum_{i=1}^m \sum_{j=1}^m B^{(i)}B^{(j)}} \int_0^h \int_0^s dW_u^{(i)} dW_s^{(j)} \|^l_p - 1 }{h},$$
%\label{multstochlog}
where the limit is taken in the sense of the existence of the generalized derivative of the Wiener process and 
$\|A\|_p$ and each $\|B^{(i)}\|_p$ are assumed to be finite.
\end{definition}
For $p=2, l\ge2$ it is easy to obtain the following estimate after proceeding as in the estimate in (\ref{msest}).
\begin{eqnarray}
\nu_2^l ( A, B^{(1:m)} )  & \le & l \mu_2 (A) + \frac{l}{2} \sum_{j=1}^m \| B^{(j)} \|_2^2  \nonumber \\
&& +\frac{l}{2} \sum_{j=1}^m ( \mu_2(B^{(j)}) + \mu_2 ( - B^{(j)} ) ) + \frac{l(l-2)}{2} \sum_{i=1}^m (\mu_2(B^{(i)}))^2
\end{eqnarray}
In general, for any $p$, it is possible to estimate $\nu_p^l ( 0, B^{(1:m)} )$ in the sense of the existence of a generalized derivative of the Wiener process such that $dW^{(i)} = \zeta^{(i)} dt,~\zeta^{(i)} \sim N(0,1),~\E(\zeta^{(i)} \zeta^{(j)}) = 0$ for $i \neq j$:
$$ \nu_p^l ( 0, B^{(1:m)} )   =   \lim_{h \to 0^+}  \frac{ \E \| I + \sum_{j=1}^m B^{(j)} \Delta W^{(j)}  + 
         \dm{\sum_{i=1}^m \sum_{k=1}^m} B^{(i)} B^{(k)} \int_0^h \int_0^s dW^{(i)}_u dW^{(k)}_s \|_p^l - 1}{h}$$ 
 $$ \le  \lim_{h \to 0^+} \bigg( \frac { \| I - (2h)\dm{\sum_{i=1}^m} \frac{1}{2} {B^{(i)}}^2 \|_p^l -1 }{(2h)}+$$ 
$$  \frac{ \E \| I +2 \Delta W \dm{\sum_{i=1}^m B^{(i)} + (\Delta W)^2  \sum_{i=1}^m {B^{(i)}}^2 + 2 \sum_{i=1}^m \sum_{j=1,~i \neq j }^m} 
B^{(i)}B^{(j)} \int_0^h \int_0^s dW_u^{(i)} dW_s^{(j)} \|^l_p -1 }{(2h)}  \bigg) $$
 $$ \le  \frac{l}{2} \mu_p(-\sum_{i=1}^m B^{(i)}) +  l \sum_{i=1}^m \| B^{(i)} \|_p  +  \frac{l}{2} \sum_{i=1}^m \| B^{(i)}  \|^2_p
 + \frac{l}{\sqrt{2}}  \sum_{i=1}^m \sum_{j=1,~i \neq j }^m \| B^{(i)}B^{(j)}\|_p.$$

Then we can upper bound the stochastic logarithmic norm as 
\begin{eqnarray}
 \nu_p^l(A,B^{(1:m)}) & \le & l \mu_p(A) -\frac{l}{2} \mu_p(\sum_{i=1}^m B^{(i)}) +  l \sum_{i=1}^m \| B^{(i)} \|_p  +  \frac{l}{2} \sum_{i=1}^m \| B^{(i)}  \|^2_p +\nonumber\\ && \frac{l}{\sqrt{2}}  \sum_{i=1}^m \sum_{j=1,~i \neq j }^m \| B^{(i)}B^{(j)}\|_p. 
\end{eqnarray} 
Similar to (\ref{nucomput1ch}) we can compute the stochastic logarithmic norm for the multi-channel case with
\begin{eqnarray}
 \nu_p^l ( A, B^{(1:m)} ) = l \E \left( \mu_p \left( A - \frac{1}{2}\sum_{i=1}^m {B^{(i)}}^2 +\sum_{i=1}^m B^{(i)} \zeta^{(i)} \right) \right), \label{directcompnu}
\end{eqnarray}
and bound it  as in (\ref{boundednu}) as
\begin{eqnarray}
 \left | \nu_p^l (A, B ^{(1:m)}\right| \le l \left \| A - \frac{1}{2} \sum_{i=1}^m {B^{(i)}}^2 \right\|_p +  l  \sum_{i=1}^m \left\| B^{(i)} \right\|_p.
\end{eqnarray}
From (\ref{directcompnu}) and as in Corollary \ref{ulnucor}, the stochastic logarithmic norm for the multiplicative multiple
channel noise can be bounded as 
\begin{eqnarray}
 l\mu_p(A) - \frac{l}{2} \sum_{i=1}^m \left( \mu_p({B^{(i)}}^2) + \mu_p(B^{(i)}) + \mu_p(-B^{(i)} )  \right) \le \nu_p^l(A, B^{(1:m)})
\nonumber \\
\le l\mu_p(A) + \frac{l}{2} \sum_{i=1}^m \left( \mu_p(-{B^{(i)}}^2) + \mu_p(B^{(i)}) + \mu_p(-B^{(i)} )\right).
\label{ulnugen}
\end{eqnarray}

\section{Relationship with Pseudospectrum}
 We have mentioned in the introduction that the logarithmic norm as a bound on the pseudospectrum of the
stability matrix 
provides an estimate of the finite time interval numerical stability of an ODE. The finite time interval numerical
stability differs from the asymptotic stability in capturing the effect of nonnormality of the stability matrices and 
local stiffness that affect the computation of the numerical integration. 
In SDEs with multiplicative noise the diffusion coefficients may 
significantly affect this transient numerical stability of an SDE. 
It may be recalled that balanced methods have been designed
\cite{milstein} to overcome difficulties arising from stiffness in both drift and diffusion. The stochastic logarithmic norm relates to the pseudospectrum of the drift coefficients
by way of diffusion coefficients acting as perturbations and thus 
captures the expected transient stability of the SDE. The estimate of stability based on the stochastic
logarithmic norm can then be used for selecting an appropriate stiff  stochastic integrator.
 \par
Let $ \| -\frac{1}{2} \sum_{i=1}^m {B^{(i)}}^2 + \sum_{i=1}^m B^{(i)} \zeta^{(i)} \|_2 = \beta $. From the definition of pseudospectrum
\cite{tref}, we may write: 
\begin{eqnarray}
\E \max_{\lambda} \Re \lambda_{\beta} (A) & = & \E \max_{\lambda} \Re{\lambda\left(A- \frac{1}{2} \sum_{i=1}^m {B^{(i)}}^2 + \sum_{i=1}^m B^{(i)} \zeta^{(i)}\right)} \nonumber \\
&& \le \E \left( \mu_2 \left(A- \frac{1}{2} \sum_{i=1}^m {B^{(i)}}^2 + \sum_{i=1}^m B^{(i)} \zeta^{(i)}\right)  \right) \nonumber \\
&&~~= \frac{1}{2}\nu_2^2 (A,B^{(1:m)})   \label{pseudodet}
\end{eqnarray}
using (\ref{directcompnu}).
For small noise with $1 \gg \beta > 0$, obviously, the stochastic logarithmic
norm gives an upper bound estimate of the mean transient numerical stability behavior 
of the deterministic ODE $dX_t = AX_t dt$.
Denoting $ \gamma =\| \sum_{i=1}^m B^{(i)} \zeta^{(i)} \|_2$, and proceeding in a similar fashion as in (\ref{pseudodet}) one obtains
\begin{eqnarray}
 \E \max_{\lambda} \Re \lambda_{\gamma} \left(A -  \frac{1}{2} \sum_{i=1}^m {B^{(i)}}^2 \right) 
  & \le & \E \left( \mu_2 \left(A- \frac{1}{2} \sum_{i=1}^m {B^{(i)}}^2 + \sum_{i=1}^m B^{(i)} \zeta^{(i)}\right) \right) \nonumber \\ && = \frac{1}{2}\nu_2^2 (A,B^{(1:m)}). \label{pseudo2} 
\end{eqnarray}
In the above inequalities the stochastic logarithmic norm appears as an upper bound on the mean maximum real part of the  perturbed spectrum of $A- \frac{1}{2} \sum_{i=1}^m {B^{(i)}}^2$ which is significant for the stochastic asymptotic stability of the SDE (\ref{mSDE}).

\section{Conclusion} This paper has extended the classical logarithmic norm to define the 
stochastic logarithmic norm for the numerical stability analysis of vector It\^o
stochastic differential equations with multi-dimensional multiplicative noise. 
Incremental estimates of the stochastic logarithmic norm due to perturbations and bounds with respect to
logarithmic norm of the drift and diffusion coefficient matrices have been studied.
Further investigation relating the stochastic logarithmic norm to the pseudo-spectrum and stiffness, both in
drift and diffusion, is needed since the stochastic logarithmic norm gives an upper bound on
the mean maximum real part of the pseudospectrum of a matrix perturbed by the noise. 
This last  property may be used in choosing stiff and balanced numerical integrators and a detailed
study for various class of integrators in this respect remains to be done.  Detailed study of application of the 
stochastic logarithmic norm to the stability analysis of non-linear SDE is necessary too.

{}


\begin{thebibliography}{99}
\bibitem{bhatia97}{\sc R.~Bhatia}, {\em Matrix Analysis}, {Springer-Verlag, 1997.}
\bibitem{burrage} {\sc K. Burrage, P. Burragea and T. Mitsui}, {\em Numerical solutions of stochastic differential equations – implementation and stability issues}, Journal of Computational and Applied Mathematics, 125:1-2 (2000), pp. 171--182 
\bibitem{burrage2} {\sc K. Burrage, P.~M.~Burrage and T.~Tian}, {\em Numerical methods for strong solutions of stochastic differential equations: an overview}, Royal Society of London Proceedings Series A, 460:2041 (2004), pp.373--402
% \bibitem{hussaini}{\sc A.~M. Croicu and M. Y. Hussaini}, 
% {\em On the expected optimal value and the optimal expected value}, 
% J. Applied Mathematics and Computation, 180(2006), pp. 330--341.
\bibitem{D59}{\sc G.~ Dahlquist}, {\em Stability and the error bounds in the numerical integration of ordinary
differential equations, Almqvist and Wiksells, Uppasala, 1958:} Transactions of Royal Institute of Technology, Stockholm, 1959.
%\bibitem{djh00}{\sc D.~J.~Higham}, {\em Mean-square and asymptotic stability of the stochastic theta method}, 
%SIAM J. Numer. Anal., 38 (2000), pp. 753--769.

\bibitem{hairer97}{\sc E. Hairer and G. Wanner}, {\em Solving Ordinary Differential Equations I, II: Stiff and 
Differential- Algebraic Problems}, {Springer-Verlag,1996.}

\bibitem{higham} {\sc D.~J. Higham and L.~N. Trefethen}, {\em Stiffness of ODEs}, BIT, 33(1993)  pp. 285.
%\bibitem{komori94}{\sc Y.~Komori, Y.~Saito and T. Mitsui}, {\em Some issues in descrite approximate solution for 
%stochastic differential equations}, Computers Math. Appl., 28(1994), pp. 269-278.
\bibitem{higham2} {\sc D.~J.~Higham and X. Mao}, {\em Nonnormality and stochastic differential equations}, BIT, 46(2006), pp. 525--532.
%\bibitem{komori95}{\sc Y.~Komori and T. Mitsui}, {\em Stable row type scheme for stochastic differential equations}, 
%Monte Carlo Methods and Appl., Computers Math. Appl., 1(1995), pp. 279-300.
\bibitem{higueras}{\sc I.~Higueras and B.~G.~Celayeta}, {\em Logarithmic norms for matrix pencils}, SIAM J. Matrix Anal., 
20 (1999), pp. 646--666.
\bibitem{kloeden}{\sc P.~E.~Kl\"oden and E.~Platen}, {\em The Numerical Solution of Stochastic Differential Equations},
Springer-Verlag, Berlin, 1992.

\bibitem{lozi58}{\sc S.~M.~Lozinskii}, {\em Error estimates for the numerical integration of ordinary differential equations}, 
part I, Izv. Vyss. Uceb. Zaved Mathematika, 6(1958), pp. 52--90.

\bibitem{mile74}{\sc G.~N.~Mil'stein},{\em Approximate integration of stochastic differential equations}, 
Theory Probab. Appl., 19(1974), pp. 557--562.

\bibitem{milstein} {\sc G.~N. Mil'stein and M.~V. Tretyakov}, {\em Stochastic Numerics for Mathematical Physics},
{Springer-Verlag, New York, 2004.}
\bibitem{raha} {\sc S. Raha and L.~R. Petzold}, {\em Constraint partitioning for stability in path-constrained dynamic optimization problems},
 SIAM J. Sci. Comput., 22 (2001), pp. 2051--2074.
\bibitem{enquist} {\sc R. Sharp, Y.-H. Tsai and B. Engquist }, {\em Multiple Time Scale Numerical Methods for the Inverted Pendulum Problem}, Lecture Notes in Computational Science and Engineering,  Springer-Verlag, 44 (2005), pp. 241-262.

\bibitem{tref} {\sc L.~N. Trefethen and M. Embree}, {\em Spectra and Pseudospectra: The Behavior of Nonnormal Matrices and Operators}, Princeton University Press, 2005.

\bibitem{mit96}{\sc Y.~Saito and T.~Mitsui}, {\em Stability analysis of numerical schemes for stochastic differential equations}, 
SIAM J. Numer. Anal., 33(1996),  pp. 2254--2267.

\bibitem{mit02}{\sc Y.~Saito and T.~Mitsui}, {\em Mean square stability of numerical schemes for stochastic differential systems}, 
Math. Sci., SIS/GSHI.,2(2002), Nagoya Univ.

\bibitem{sanzserna} {\sc J.~M. Sanz-Serna}, {\em Modulated Fourier expansions and heterogeneous multiscale methods}, 
IMA Journal of Numerical Analysis, doi:10.1093/imanum/drn031 (2008)

\bibitem{strom75}{\sc T.~ Str\^{o}m}, {\em On logarithmic norms}, SIAM J. Numer. Anal., 12(1975), pp. 741--753.

\bibitem{gustaf06}{\sc G.~ S\"oderlind}, {\em The logarithmic norm history and modern theory}, BIT, 46 (2006), pp. 631--652.

\end{thebibliography}
\end{document}